\documentclass[12pt,oneside,english]{amsart}
\usepackage[latin1]{inputenc}
\usepackage{amssymb}

\makeatletter
\newtheorem{theorem}{Theorem}
\newtheorem{lemma}{Lemma}

\newtheorem{definition}{Definition}
\newtheorem{condition}{Condition}
\newtheorem{remark}{Remark}

\newtheorem{problem}{Problem}
\newtheorem{question}{Question}

\usepackage{babel}

\makeatother
\begin{document}
\baselineskip=17pt

\title[On small intervals containing primes]{On small intervals containing primes}

\author{Vladimir Shevelev}
\address{Departments of Mathematics \\Ben-Gurion University of the
 Negev\\Beer-Sheva 84105, Israel. e-mail:shevelev@bgu.ac.il}


\begin{abstract}
 Let $p$ be an odd prime, such that $p_n<p/2<p_{n+1},$ where $p_n$ is the n-th prime. We study the following question: with what probability does there exist a prime in the interval $(p, 2p_{n+1})?$  After the strong  definition of the probability with help of the Ramanujan primes ([11], [12])and the introducing pseudo-Ramanujan primes, we show, that if such probability $\mathbf{P}$ exists, then $\mathbf{P}\geq0.5.$ We also study a symmetrical case of the left intervals, which connected with sequence A080359 in [10].
 \end{abstract}

\maketitle

\section{Introduction}
As well known, the Bertrand's postulate (1845) states that, for $x>1,$ always there exists a prime in interval $(x, 2x).$ This postulate very quickly-five years later- became a theorem due to Russian mathematician P.L.Chebyshev (cf., e.g., [9, Theorem 9.2]). In 1930 Hoheisel[3] proved that, for $x>x_0(\varepsilon),$ the interval $(x, x+x^{1-\frac {1} {33000}+\varepsilon}]$ always contains a prime. After that there were a large chain of improvements of the Hoheisel's result. Up to now, probably, the best known result belongs to Baker, Harman abd Pintz[1], who showed that even the interval $(x, x+x^{0.525})$ contains a prime. Their result is rather close to the best result which gives the Riemann
hypothesis: $p_{n+1}-p_n=O(\sqrt{p_n}\ln p_n)$ (cf. [4, p.299]), but still very far from the Cram\'{e}r's 1937 conjecture which states that already the interval $(x,x+(1+\varepsilon)\ln^2x]$ contains a prime for sufficiently large $x.$\newline
\indent Everywhere during this paper we understand that $p_n$ is the n-th prime.\newline
Let $p$ be an odd prime. Let, furthermore, $p_n<p/2<p_{n+1}.$ According to the Bertrand's postulate, between $p/2$ and $p$ there exists a prime. Therefore, $p_{n+1}\leq p.$ Again, by the Bertrand's postulate, between $p$ and $2p$ there exists a prime. More subtle question is the following.
\begin{question}
 Let $p$ be a randomly chosen odd prime. Suppose that $p/2$ lies in the interval $(p_n,\enskip p_{n+1}).$ With what probability does there exist a prime in the interval $(p, 2p_{n+1})?$
   \end{question}
   At the first we should formulate more exactly what we understand under such probability. To this end we start with
    two conditions for odd primes and their equivalence. An important role for our definition of the desired probability play Ramanujan primes ([11]-[12]) and also Pseudo-Ramanujan primes which we introduce below.
    \section{Equivalence of two conditions for odd primes}
Consider the following two conditions for primes:
\begin{condition}\label{1}
Let $p=p_n,$ with  $n>1.$  Then all integers $(p+1)/2, (p+3)/2, ... , (p_{n+1}-1))/2 $ are composite numbers.
\end{condition}
\begin{condition}\label{2}
Let, for an odd prime $p,$ we have $p_m< p/2< p_{m+1}.$ Then the interval  $(p,\enskip 2p_{m+1})$ contains a prime.
\end{condition}
\begin{lemma}\label{1}
 Conditions 1 and 2 are equivalent.
 \end{lemma}
 \bfseries Proof. \mdseries If Condition 1 is valid, then $p_{m+1}>(p_{n+1}-1)/2,$ i.e. $p_{m+1}\geq(p_{n+1}+1)/2.$ Thus $2p_{m+1}> p_{n+1}>p_n=p,$ and Condition 2 is valid; conversely, if Condition 2 satisfies, i.e. $p_{m+1}>p/2$
 and   $2p_{m+1}>p_{n+1}>p=p_n.$ If $k$ is the least positive integer, such that $p_m<p_n/2<(p_n+k)/2<(p_{n+1}-1)/2$ and $(p_n+k)/2$ is prime, then $p_{m+1}=(p_n+k)/2$ and $p_{n+1}-1>p_n+k=2p_{m+1}>p_{n+1}.$ Contradiction shows that
 Condition 1 is valid. $\blacksquare$
 \section{ Ramanujan primes}
  In 1919 S. Ramanujan [7]-[8] unexpectedly gave a new short and elegant proof of the Bertrand's postulate. In his proof appeared a sequence of primes
\begin{equation}\label{1}
2,11,17,29,41,47,59, 67, 71, 97, 101, 107, 127, 149, 151, 167,...
 \end{equation}
 For a long time, this important sequence was not presented in the Sloane's OEIS [9]. Only in 2005 J. Sondow published it in OEIS (sequence A104272).
\begin{definition}\label{1} (J. Sondow$[10]$)For $n \geq 1$, the \upshape nth Ramanujan  prime \slshape is the smallest positive integer  $(R_n)$  with the property that if $x\geq R_n,$  then $\pi(x) -\pi(x/2)\geq n.$
  \end{definition}
In [11], J. Sondow obtained some estimates for $R_n$ and, in particular,  proved that, for every $n>1,\enskip R_n>p_{2n}.$ Further, he proved that for $n\rightarrow\infty, \enskip R_n\sim p_{2n}.$  From this, denoting $\pi_R(x)$ the counting function of the Ramanujan  primes not exceeding $x,$ we have $R_{\pi_R(x)}\sim2\pi_R(x)\ln \pi_R(x).$
Since $R_{\pi_R(x)}\leq x<R_{\pi_{R(x)}+1},$
then $ x\sim p_{2\pi_R(x)}\sim2\pi_R(x)\ln \pi_R(x),$ as $x\rightarrow\infty,$ and we conclude that\newpage
 \begin{equation}\label{2}
\pi_R(x)\sim \frac {x} {2\ln x}.
 \end{equation}
 It is interesting that quite recently S. Laishram (see [10], comments to A104272) has proved a Sondow conjectural inequality  $ R_n < p_{3n}$ for every positive $n.$
\section{ Satisfaction Conditions 1,2 for Ramanujan primes}
 \begin{lemma}\label{2}
 If p is an odd Ramanujan prime, then Conditions 1 and 2 satisfy.
 \end{lemma}
 \bfseries Proof. \mdseries In view of Lemma 1, it is sufficient to prove that Condition 1 satisfies. If Condition 1
 does not satisfy, then suppose that $p_m=R_n<p_{m+1}$ and  $k$ is the least positive integer, such that $q=(p_m+k)/2$ is prime not more than $(p_{m+1}-1)/2.$ Thus
  \begin{equation}\label{3}
  R_n=p_m<2q<p_{m+1}-1.
 \end{equation}
 As Sondow proved ([12]), $R_n-1$ is the maximal integer for which the equality
  \begin{equation}\label{4}
 \pi(R_n-1)-\pi((R_n-1)/2)=n-1
 \end{equation}
 holds. However, according to (5), $\pi(2q)= \pi(R_n-1)+1$ and in view if the minimality of the prime $q,$ in the interval $((R_n-1)/2,q)$ there are not any prime. Thus $\pi(q)=\pi((R_n-1)/2)+1$ and
 $$ \pi(2q)-\pi(q)=\pi(R_n-1)-\pi((R_n-1)/2)=n-1.$$
 Since, by (5), $2q>R_n,$ then this contradicts to the property of the maximality of $R_n$ in (6). $ \blacksquare$\newline\newline
 \indent Note that, there are non-Ramanujan primes which satisfy Conditions 1,2. We call them  \slshape pseudo-Ramanujan \upshape primes $(PR)_n.$ The first terms of the sequence of pseudo-Ramanujan primes are:
\begin{equation}\label{5}
109,137,191,197,283,521,...
 \end{equation}
 \begin{definition}\label{2}
 We call a prime $p$ an \upshape RPR-prime \slshape if $p$ satisfies Condition 1 (or, equivalently, Condition 2).
  \end{definition}
  From the above it follows that a RPR-prime is either Ramanujan or pseudo-Ramanujan prime.
  Thus we see that the relative density (if it exists) of RPR-primes (and only of them) with respect to all primes not
  exceeding $N,$ for $N$ tends to the infinity should give the answer on Question 1. More\newpage exactly, denote $(RPR)_n$ the $n$-th pseudo-Ramanujan prime and $\pi_{RPR}(n)$ the number of $RPR$-primes not exceeding $n.$
\begin{definition}\label{3}
 Let $p_n<p/2<p_{n+1}.$ Under the probability $\mathbf{P}$ that there exists a prime in the interval $(p, 2p_{n+1})$ we understand, if it exists, the limit
 $$\mathbf{P}:=\lim_{n\rightarrow\infty} \frac  {\pi_{RPR}(n)} {\pi (n)},$$ or, the same, by Prime Number Theorem,
 $$\mathbf{P}=\lim_{n\rightarrow\infty} \frac  {\pi_{RPR}(n)} {{n}/\ln n}.$$

  \end{definition}
  \section{A sieve for selection RPR-primes from all primes}
Denote $\mathbb{PR}$(respectively, $\mathbb{RPR})$ the set of all pseudo-Ramanujan primes (respectively, RPR-primes).
The probability under consideration is
$$\mathbf{P}(x\in\mathbb{R}/x\in\mathbb{P})+\mathbf{P}(x\in\mathbb{PR}/x\in\mathbb{P})=
\mathbf{P}(x\in\mathbb{RPR}/x\in\mathbb{P}).$$
Therefore, it is interesting to build a sieve for selection RPR-primes from all primes. Recall that the Bertrand sequence $\{b(n)\}$ is defined as $ b(1)=2,$ and, for $n\geq2,\enskip b(n)$ is the largest prime less than $2b(n-1)$ (see A006992 in [10]):
\begin{equation}\label{6}
2, 3, 5, 7, 13, 23, 43,...
 \end{equation}
Put
\begin{equation}\label{7}
B_1=\{b^{(1)}(n)\}=\{{b(n)}\}.
 \end{equation}
Further we build sequences $B_2=\{b^{(2)}(n)\}, B_3=\{b^{(3)}(n)\},...$ according the following inductive rule: if we have sequences $B_1,...,B_{k-1},$ let us consider the minimal prime $p^{(k)}\not \in \bigcup_{i=1}^{k-1}B_i.$ Then the sequence $\{b^{(k)}(n)\}$ is defined as $b^{(k)}(1)=p^{(k)},$ and, for $n\geq2,\enskip b^{(k)}(n) $ is the largest prime less than $2b^{(k)}(n-1).$ So, we obtain consequently:
\begin{equation}\label{8}
B_2=\{11, 19, 37, 73,...\}
 \end{equation}
 \begin{equation}\label{9}
B_3=\{17, 31, 61, 113,...\}
 \end{equation}
 \begin{equation}\label{10}
B_3=\{29, 53, 103, 199,...\}
 \end{equation}
 etc., such that, putting $p^{(1)}=2,$ we obtain the sequence
 \begin{equation}\label{11}
\{p^{(k)}\}_{k\geq1}=\{2, 11, 17, 29, 41, 47, 59, 67, 71, 97, 101, 107, 109, 127, ...\}
 \end{equation}
 Sequence (11) coincides with sequence (1) of the Ramanujan primes up to the $12$-th term, but the $13$-th term of this sequence is $109 $ which is the first term of sequence (5) of the pseudo-Ramanujan primes.\newpage
 \begin{theorem}\label{1} For $n\geq1,$ we have
 \begin{equation}\label{12}
p^{(n)}=(RPR)_n
 \end{equation}
 where $(RPR)_n$ is the $n$-th RPR-prime.
\end{theorem}
\bfseries Proof. \mdseries The least omitted prime in (7) is $p^{(2)}=11=(RPR)_2$;  the least omitted prime in the union of (7) and (8) is $p^{(3)}=17=(RPR)_3.$ We use the induction. Let we have already built primes $$p^{(1)}=2, p^{(3)},...,p^{(n-1)}=(RPR)_{n-1}.$$ Let $q$ be the least prime which is omitted in the union $\bigcup_{i=1}^{n-1}B_i,$  such that $q/2$ is in interval $(p_{m}, p_{m+1}).$ According to our algorithm, $q$ which is dropped should not be the large prime in the interval $(p_{m+1}, 2p_{m+1}).$ Then there are primes in the interval  $q, 2p_{m+1})$; let $r$ be one of them. Then we have $2p_m<q<r<2p_{m+1}.$ This means that $q,$ in view of its minimality between the dropping primes more than $(RPR)_{n-1}=p^{(n-1)}$, is the least $RPR$-pime more than $(RPR)_{n-1}$ and  the least prime of the form $p^{(k)}$ more than $p^{(n-1)}$. Therefore, $q=p^{(n)}=(RPR)_{n}.\blacksquare$\newline\newline
\indent Unfortunately the research of this sieve seems much more difficult than the research of the Eratosthenes one for primes. For example, the following question remains open.
\begin{problem}
 With help of the sieve of Theorem 1 to find a formula for the counting function of RPR-primes not exceeding $x.$
   \end{problem}
 The following theorem is proved even without the supposition of the existing the limit in Definition 3.
 \begin{theorem}\label{2} Denote $$\underline{\mathbf{P}}=\lim\inf_{n\rightarrow\infty}\frac {\pi_{RPR}(n)}{\pi(n)}$$ ("lower probability"). Then we have
$$\underline{\mathbf{P}}\geq\frac {1} {2}.$$

\end{theorem}
\bfseries Proof. \mdseries Using (2), we have
$$ \underline{\mathbf{P}}= {\lim\inf}_{n\rightarrow\infty }\pi_{RPR}(n)/\pi(n)\geq{\lim}_{n\rightarrow\infty }\pi_{R}(n)/\pi(n)=1/2. \blacksquare$$

\indent  D. Berend [2] gave another very elegant proof of this theorem.

\bfseries Second proof of Theorem 2. \mdseries We saw that if the interval  $(2p_m, 2p_{m+1})$ with odd $p_m$ contains a prime $p,$ then the interval $(p, 2p_{m+1})$ contains in turn a prime if and only if $p$ is a RPR-primes. Let $n\geq7.$ In the range from 7 up to $n$ there are $\pi(n)-3$ primes. Put\newpage
\begin{equation}\label{13}
 h=h(n)=\lfloor(\pi(n)-1)/2\rfloor.
 \end{equation}\label{13}
   Look at $h$ intervals:
\begin{equation}\label{14}
(2p_2, 2p_3), (2p_3, 2p_4), ..., (2p_{h+1}, 2p_{h+2}).
\end{equation}
Our $\pi(n)-3$ primes are somehow distributed in these $h$
intervals. Suppose $k=k(n)$ of the intervals contain at least one prime and
$h-k$ contain no primes. Then for exactly $k$ primes there is no primes
between them and the next $2p_j,$ and for the other $\pi(n)-3-k $ there is. Hence,
among $\pi(n)-3$ primes exactly $\pi(n)-3-k$ are RPR-primes and exactly $k$ non-RPR-primes. Therefore, since $k(n)\leq h(n)\leq \pi(n)/2,$ then for the desired lower probability that there is a prime we have:
\begin{equation}\label{15}
\underline{\mathbf{P}}={\lim\inf}_{n\rightarrow\infty} {\frac  {\pi_{RPR}(n)} {\pi(n)-3}}={\lim\inf}_{n\rightarrow\infty}{\frac {\pi(n)-k(n)} {\pi(n)}}\geq1/2.
\end{equation}
$\blacksquare$\newline
\section{A heuristic idea of Greg Martin }
 Greg Martin [5] proposed the following heuristic arguments, which show that $\mathbf{P}$ is close to $2/3.$
  "Imagine the following process: start from $p $ and examine the
 numbers $p+1, p+2, ... $ in turn.  If the number we're examining is
 odd, check if it's a prime: if so, we "win".  If the number we're
 examining is twice an odd number (that is,$ 2 \pmod4$), check if it's twice a prime:
 if so, we "lose".  In this way we "win" if and only if there
 is a prime in the interval $(p, 2p_{n+1})$, since we either find such a prime
 when we "win" or else detect the endpoint $2p_{n+1},$ when we "lose".

 Now if the primes were distributed totally randomly, then
 the probability of each odd number being prime would be the same(roughly $1/ln p$), while the probability of a $2 \pmod 4$ number being twice a prime would be roughly $1/ln(p/2)$, which for p large is about the
 same as $1/ln p.$  However, in every block of 4 consecutive
 integers, we have
 two odd numbers that might be prime and only one $2 \pmod 4$
 number that might be twice a prime. Therefore we expect that we "win"
 twice as often as we "lose", since the placement of primes should
 behave statistically randomly in the limit; in other words, we expect
 to "win" $\mathbf{P}_0=2/3$ of the time."
 His computions what happens for $p$ among the first million primes show that
 the probability of "we win" has a steadily increasing trend as
 $p$ increases, and among the first million primes about $61.2 $ of
 them have a prime in the interval $(p, 2p_{n+1}).$\newline\newline

\begin{remark} The following heuristic arguments are seductive but wrong:
\end{remark}
Consider the probability that a random interval $I$ from system (14) contains a prime. In order to say about a statistics, consider $I$ with "average" length, which for large $n$ equals to $2p_n/(n-1)\sim 2\ln n.$ Note that, the probability that a random integer from $(0, \enskip 2p_{n+1})$ is prime is $1/\ln n.$ Thus the proportion  $\theta$ of the absence a prime in $I$ is \enskip $\theta=(1-\frac {1} {\ln n})^{2\ln n}=e^{-2}(1+o(1)).$
Therefore, we expect that the number of intervals (14) containing a prime is approximately $(1-\theta)n.$ Now using (15) for $k=\theta n,$ we obtain the probability $$\mathbf{P}_1=1/2(1+e^{-2})=0.5676676...\enskip.$$
It is  \bfseries not correct \mdseries since, as well known (see [6, Chapter 5]), for $n$ tends to the infinity and any $c>0,$ there is a finite part of differences of the consecutive primes $p_{n+1}-p_n$ which are less than $c\ln n.$ This makes $\theta$ less than $e^{-2}$ and  $\mathbf{P}_1$ more than $1/2(1+e^{-2}).$

\section{A symmetrical case of the left intervals} It is clear that for the symmetrical problem of the existence a prime in the left interval $(2p_n, p)$ (for the same condition $p_n<p/2<p_{n+1}$) we have the same results.
Therefore, this case is not interesting from the formal-probabilistic point of view, but it is more interesting from the sequences point of view. Indeed, now in our construction the role of the Ramanujan primes play other primes which appear in OEIS [9] earlier (2003) than the Ramanujan primes due to E. Labos (see sequence A080359):
\begin{equation}\label{18}
2, 3, 13, 19, 31, 43, 53, 61, 71, 73, 101, 103, 109, 113, 139, 157, 173,...
 \end{equation}
 These primes we call the \slshape Labos primes. \upshape
\begin{definition}\label{5}(cf. [$9$, A080359] For $n \geq 1$, the \upshape nth Labos  prime\slshape \enskip is the smallest positive integer  $(L_n)$  for which $\pi(L_n)-\pi(L_n/2)=n.$
  \end{definition}
  Note that, since ([11])
  \begin{equation}\label{19}
\pi(R_n)-\pi(R_n/2)=n,
 \end{equation}
 then, by the Definition 2, we have
\begin{equation}\label{20}
L_n\leq R_n.
 \end{equation}
 As above, one can prove the equivalence of the following conditions on primes:

\begin{condition}\label{3}
Let $p=p_n$ with $n\geq3.$ Then all integers $(p-1)/2, (p-3)/2, ... , (p_{n-1}+1)/2 $ are composite numbers.
\end{condition}\newpage
\begin{condition}\label{4}
Let $p_m< p/2< p_{m+1}.$ Then the interval  $(2p_m,\enskip p)$ contains a prime.
\end{condition}
Furthermore, by the same way as for Lemma 2, one can prove that if $p$ is a Labos  prime, then Conditions 3 and 4 satisfy. But again there are non-Labos primes which satisfy Conditions 3,4. We call them  \slshape pseudo-Labos\upshape \enskip primes $(PR)_n.$ The first terms of the sequence of pseudo-Labos primes are:
\begin{equation}\label{21}
131, 151, 229, 233, 311, 571 ,...
 \end{equation}
 \begin{definition}\label{6}
 We call a prime $p$ a \upshape LPL-prime \slshape if $p$ satisfies Condition 3 (or, equivalently, Condition 4).
  \end{definition}
  From the above it follows that a LPL-prime is either Labos or pseudo-Labos prime.
  Note that for the LPL-primes one can build a sieve with help of the Sloan's primes (see A055496 [10]) and the corresponding generalizations of them (cf. constructing in Section 5).

\bfseries Acknowledgment.\mdseries \enskip The author is grateful to Daniel Berend (Ben Gurion University, Israel) for important private communication, Greg Martin (University of British Columbia, Canada) and Jonathan Sondow (USA) for their helpful comments which promoted a significant improvement of the paper.

\end{document}